\documentclass[12pt]{amsart} 
\usepackage[active]{srcltx}
\usepackage{color}
\usepackage[english]{babel}
\usepackage[latin1]{inputenc}
\usepackage{amssymb}
\usepackage{amsmath}
\usepackage{amscd}
\usepackage{latexsym}
\usepackage{amsthm}
\usepackage{mathrsfs}

\theoremstyle{plain}
\newtheorem{thm}{Theorem}[section]
\newtheorem{cor}[thm]{Corollary}
\newtheorem{lem}[thm]{Lemma}
\newtheorem{prop}[thm]{Proposition}
\newtheorem{conj}[thm]{Conjecture}
\newtheorem{defn}[thm]{Definition}
\newtheorem{oss}[thm]{Remark}

\def\vol{\mbox{vol}}

\begin{document}

\title[K\"ahlerness of moduli spaces]{K\"ahlerness of moduli spaces of stable sheaves over non-projective K3 surfaces}
\author[Perego]{Arvid Perego}

\begin{abstract}
We show that a moduli space of slope-stable sheaves over a K3 surface is an irreducible hyperk\"ahler manifold if and only if its second Betti number is the sum of its Hodge numbers $h^{2,0}$, $h^{1,1}$ and $h^{0,2}$.
\end{abstract}

\date{\today}
\thanks{
 }

\maketitle

\tableofcontents

\clearpage

\section{Introduction}

Irreducible hyperk\"ahler manifolds are compact, connected K\"ahler manifolds which are simply connected, holomorphically symplectic and have $h^{2,0}=1$. Very few examples of them are know up to today, and all the known deformation classes arise from moduli spaces of semistable sheaves on a projective K3 surface or on an abelian surface.

In \cite{PT} we showed that if $S$ is any K3 surface, the moduli space $M_{v}^{\mu}(S,\omega)$ of $\mu_{\omega}-$stable sheaves on $S$ of Mukai vector $v=(r,\xi,a)\in H^{2*}(S,\mathbb{Z})$ is a compact, connected complex manifold, it carries a holomorphic symplectic form and it is of K3$^{[n]}-$type (i. e. it is deformation equivalent to a Hilbert scheme of points on a projective K3 surface). This holds under some hypothesis on $\omega$ and $v$ (the K\"ahler class $\omega$ has to be $v-$generic, and $r$ and $\xi$ have to be prime to each other: we refer the reader to \cite{PT} for the definition of $v-$genericity).

The main open question about these moduli spaces is if they carry a K\"ahler metric: if it is so, it follows that they are all irreducible hyperk\"ahler manifolds of K3$^{[n]}-$type. The answer to this question is affirmative at least in three cases: when $S$ is projective; when $M^{\mu}_{v}(S,\omega)$ is a surface; when $M^{\mu}_{v}(S,\omega)$ parametrizes only locally free sheaves. This lead us to the following

\begin{conj}
\label{conj:con}
The moduli spaces $M_{v}^{\mu}(S,\omega)$ are K\"ahler manifolds.
\end{conj}

Evidences are provided by the previous examples where the moduli spaces are indeed K\"ahler, and by the fact that their geometry is somehow similar to that of an irreducible hyperk\"ahler manifold: in \cite{PT} we show that on their second integral cohomology there is a non-degenerate quadratic form defined as the Beauville form of irreducible hyperk\"ahler manifolds.

But still, this analogy is not sufficient to guarantee that the moduli spaces are K\"ahler: it is known since \cite{Gua1}, \cite{Gua2} and \cite{Gua3} that there are examples of compact, simply connected, holomorphically symplectic manifolds having $h^{2,0}=1$ which are not K\"ahler, but their second integral cohomology carries a non-degenerate quadratic form, and the Local Torelli Theorem holds.

The aim of this paper is to show that the previous conjecture holds true under some additional hypothesis on the second Betti number of $M^{\mu}_{v}(S,\omega)$.

\subsection*{Acknowledgements}
The author wishes to thank M. Toma for several useful conversations during the preparation of this work.

\subsection{Main definitions and notations}

In this section we collect all the definitions and notations we will use in the following.

\begin{defn}
A \textbf{holomorphically symplectic manifold} is a complex manifold which carries an everywhere non-degenerate holomorphic closed $2-$form (called \textbf{holomorphic symplectic form}).
\end{defn}

We notice that a compact holomorphically symplectic manifold is always of even complex dimension, and a holomorphic form $\sigma$ defines an isomorphism $\sigma:T_{X}\longrightarrow\Omega_{X}$ of vector bundles, where $T_{X}$ is the tangent bundle of $X$, and $\Omega_{X}$ is the cotangent bundle of $X$ (i. e. the dual bundle of $T_{X}$).

Let $X$ be a compact, connected complex manifold of complex dimension $d$, and $k\in\{0,...,2d\}$ and $p,q\in\mathbb{N}$.

\begin{defn}
The $k-$\textbf{Betti number} of $X$ is $$b_{k}(X)=dim_{\mathbb{C}}H^{k}(X,\mathbb{C}),$$and the $(p,q)-$\textbf{th Hodge number} of $X$ is $$h^{p,q}(X)=dim_{\mathbb{C}}H^{q}(X,\Omega_{X}^{p}).$$
\end{defn}

The general relation between the Betti and the Hodge numbers of $X$ is that $$b_{k}(X)\leq\sum_{p+q=k}h^{p,q}(X)$$for every $k\in\{0,...,2d\}$, and the equality holds for every $k$ if and only if the Fr\"olicher spectral sequence of $X$ degenerates at the $E_{1}$ level.

\begin{defn}
A complex manifold $X$ is a $b_{2}-$\textbf{manifold} if $$b_{2}(X)=h^{2,0}(X)+h^{1,1}(X)+h^{0,2}(X).$$
\end{defn}

\begin{defn}
A compact, connected complex manifold $X$ \textbf{is in the Fujiki class $\mathcal{C}$} if $X$ is bimeromorphic to a compact K\"ahler manifold.
\end{defn} 

Among all manifolds in the Fujiki class $\mathcal{C}$ we clearly have compact K\"ahler manifolds. Moreover, a compact, connected complex manifold in the Fujiki class $\mathcal{C}$ verifies the $\partial\overline{\partial}-$lemma, and hence the Fr\"ohlicher spectral sequence of $X$ degenerates at the $E_{1}$ level. In particular, it is a $b_{2}-$manifold. 

If $f:\mathcal{X}\longrightarrow B$ is a holomorphic fibration, for every $b\in B$ we let $X_{b}:=f^{-1}(b)$. Let $X$ be a compact, connected complex manifold and $B$ a connected complex manifold.

\begin{defn}
A \textbf{deformation of} $X$ \textbf{along} $B$ is a smooth and proper family $f:\mathcal{X}\longrightarrow B$ such that there is $0\in B$ such that $X_{0}$ is biholomorphic to $X$.
\end{defn}

Let now $\mathcal{P}$ be a property of complex manifolds. 

\begin{defn}
We say that the property $\mathcal{P}$ is \textbf{open in the analytic topology} (resp. \textbf{in the Zariski topology}) if for every deformation $\mathcal{X}$ along a connected complex manifold $B$, the set of those $b\in B$ such that $X_{b}$ verifies $\mathcal{P}$ is an analytic (resp. Zariski) open subset of $B$;
\end{defn}

K\"ahlerness is open in the analytic topology, but in general it is neither open in the Zariski topology, nor closed. Being a $b_{2}-$manifold is open in the Zariski topology. Being in the Fujiki class $\mathcal{C}$ is not open in general.

The last definitions we need are the following:

\begin{defn}
Let $X$ be an compact, connected complex manifold.
\begin{enumerate}
 \item The manifold $X$ is \textbf{irreducible hyperk\"ahler} if it is a K\"ahler manifold which is simply connected, holomorphically simplectic and $h^{2,0}(X)=1$. 
 \item The manifold $X$ is \textbf{deformation equivalent to an irreducible hyperk\"ahler manifold} if there is a connected complex manifold $B$ and a deformation $\mathcal{X}\longrightarrow B$ of $X$ along $B$ for which there is $b\in B$ such that $X_{b}$ is an irreducible hyperk\"ahler manifold.
 \item The manifold $X$ is \textbf{limit of irreducible hyperk\"ahler manifolds} if there is a smooth and proper family $\mathcal{X}\longrightarrow B$ along a smooth connected base $B$ such that and a sequence $\{b_{n}\}$ of points of $B$ converging to $0$, such that $X_{b_{n}}$ is an irreducible hyperk\"ahler manifold.
\end{enumerate}
\end{defn}

\subsection{Main results and structure of the paper}

The main result of the paper is the following

\begin{thm}
\label{thm:main}Let $S$ be a K3 surface, $\omega$ a K\"ahler class on $S$. We let $v=(r,\xi,a)\in H^{2*}(S,\mathbb{Z}$) be such that $r>0$ and $\xi\in NS(S)$. Suppose that $r$ and $\xi$ are prime to each other, and that $\omega$ is $v-$generic. Then the moduli space $M=M^{\mu}_{v}(S,\omega)$ of $\mu_{\omega}-$stable sheaves on $S$ with Mukai vector $v$ is K\"ahler if and only if it is a $b_{2}-$manifold.
\end{thm}

Clearly, this proves Conjecture \ref{conj:con} under the additional hypothesis of the moduli spaces being $b_{2}-$manifolds. This has an immediate corollary:

\begin{cor}
\label{cor:kahlerzariski}Let $\mathcal{M}\longrightarrow B$ be any smooth and proper family of moduli spaces of sheaves verifying the conditions of Theorem \ref{thm:main}. The set of $b\in B$ such that $M_{b}$ is K\"ahler is a Zariski open subset of $B$.
\end{cor}

In view of of Theorems 1.1 and 1.2 of \cite{PT}, another immediate corollary is the following:

\begin{cor}
\label{cor:hyperk}Let $S$ be a K3 surface, $\omega$ a K\"ahler class on $S$ and $v=(r,\xi,a)\in H^{2*}(S,\mathbb{Z}$) be such that $r>0$ and $\xi\in NS(S)$. Suppose that $r$ and $\xi$ are prime to each other, that $\omega$ is $v-$generic, and that $M^{\mu}_{v}(S,\omega)$ is a $b_{2}-$manifold.
\begin{enumerate}
 \item The moduli space $M^{\mu}_{v}(S,\omega)$ is an irreducible hyperk\"ahler manifold of K3$^{[n]}-$type, which is projective if and only if $S$ is projective.
 \item If $v^{2}\geq 2$, there is a Hodge isometry $\lambda_{v}:v^{\perp}\longrightarrow H^{2}(M^{\mu}_{v},\mathbb{Z})$.
\end{enumerate}
\end{cor}

The case $v^{2}=0$ was already treated in \cite{PT}: in this case there is a Hodge isometry $$\lambda_{v}:v^{\perp}/\mathbb{Z}v\longrightarrow H^{2}(M^{\mu}_{v},\mathbb{Z}),$$and there is no need to suppose that $M^{\mu}_{v}(S,\omega)$ is a $b_{2}-$manifold here. For as a consequence Theorem 1.1 of \cite{PT} we already know that $M^{\mu}_{v}(S,\omega)$ is a K3 surface.

\par\bigskip

The proof of Theorem \ref{thm:main} is an application of general results about compact, connected complex $b_{2}-$manifolds which are holomorphically symplectic and limit of irreducible hyperk\"ahler manifolds. The starting point is the following:

\begin{thm}
\label{thm:thm1}Let $X$ be a compact, connected holomorphically symplectic $b_{2}-$manifold which is deformation equivalent to an irreducible hyperk\"ahler manifold. Then on $H^{2}(X,\mathbb{Z})$ there is a non-degenerate quadratic form $q_{X}$ of signature $(3,b_{2}(X)-3)$, and the Local Torelli Theorem holds.
\end{thm}

This result is due to Guan \cite{Gua3}, and it is a generalization of the well-known analogue for irreducible hyperk\"ahler manifolds proved by Beauville in \cite{Beau}. By Local Torelli Theorem we mean that the period map is locally a biholomorphism (as in the case of irreducible hyperk\"ahler manifolds). We will recall the definition of the Beauville form and the Local Torelli Theorem in Section 2.

\begin{oss}
\label{oss:localtorelli}
{\rm In \cite{PT} we proved (see Theorem 1.1 there) that if $M$ is a moduli spaces of slope-stable sheaves over a non-projective K3 surface (veryfing all the hypothesis of Theorem \ref{thm:main}), then on $H^{2}(M,\mathbb{Z})$ there is a non-degenerate quadratic form of signature $(3,b_{2}(M)-3)$. This is proved without any assumption on $b_{2}(M)$. Anyway, in \cite{PT} we have not proved the Local Torelli Theorem, and here we are able to prove it only by assuming $M$ to be a $b_{2}-$manifold.}
\end{oss}

In Section 3 we consider compact, connected holomorphically symplectic $b_{2}-$manifolds which are not only deformation equivalent to an irreducible symplectic manifold, but which are moreover limit of irreducible hyperk\"ahler manifolds. The main result of Section 3 is the following:

\begin{thm}
\label{thm:thm2}Let $X$ be a compact, connected holomorphically symplectic $b_{2}-$manifold which is limit of irreducible hyperk\"ahler manifolds. Then $X$ is bimeromorphic to an irreducible hyperk\"ahler manifold (hence, it is in the Fujiki class $\mathcal{C}$).
\end{thm}

As we will see, a moduli space $M$ verifying the hypothesis of Theorem \ref{thm:main} is a compact, connected holomorphically symplectic $b_{2}-$manifold which is limit of irreducible hyperk\"ahler manifolds. As a consequence, asking for these moduli spaces to be $b_{2}-$manifolds is enough to conclude that they are all bimeromorphic to a compact irreducible hyperk\"ahler manifold.

The proof of Theorem \ref{thm:thm2} is based on a well-known strategy already used by Siu in \cite{Siu} to show that all K3 surfaces are K\"alher, and by Huybrechts in \cite{Hu99} to show that non-separated marked irreducible hyperk\"ahler manifolds are in fact bimeromorphic. More precisely, if $\Lambda$ is a lattice, we say that a compact complex manifold $X$ \textit{carries a $\Lambda-$marking} if on $H^{2}(X,\mathbb{Z})$ there is non-degenerate quadratic form, and there is an isometry $\phi:H^{2}(X,\mathbb{Z})\longrightarrow\Lambda$. The pair $(X,\phi)$ is called a \textit{$\Lambda-$marked manifold}.

The set of (equivalence classes of) $\Lambda-$marked manifolds is denoted $\mathcal{M}_{\Lambda}$: as a consequence of Theorem \ref{thm:thm1}, it contains the subset $\mathcal{M}_{\Lambda}^{s}$ of $\Lambda-$marked manifolds $(X,\phi)$ where $X$ is a compact holomorphically symplectic $b_{2}-$manifold which is deformation equivalent to an irreducible hyperk\"ahler manifold (and whose Beauville lattice is isometric to $\Lambda$).

By the Local Torelli Theorem we can give $\mathcal{M}_{\Lambda}^{s}$ the structure of complex space, in which we have a (non-empty) open subset $\mathcal{M}_{\Lambda}^{hk}$ of irreducible hyperk\"ahler manifolds. We let $\overline{\mathcal{M}}_{\Lambda}^{hk}$ be its closure in $\mathcal{M}_{\Lambda}^{s}$.

Theorem \ref{thm:thm2} can be restated by saying that if $(X,\phi)\in\overline{\mathcal{M}}_{\Lambda}^{hk}$, then $X$ is bimeromorphic to an irreducible hyperk\"ahler manifold. This is the statement we prove: the idea of the proof is that if $(X,\phi)\in\overline{\mathcal{M}}_{\Lambda}^{hk}$, then $(X,\phi)$ is non-separated from a point $(Y,\psi)\in\mathcal{M}^{hk}_{\Lambda}$. A standard argument shows that $X$ and $Y$ have to be bimeromorphic.

Theorem \ref{thm:thm2} is just an intermediate result on the way to the K\"ahlerness of the moduli spaces, and it is used in Section 4 to prove that on a compact, connected holomorphically symplectic $b_{2}-$manifold which is limit of irreducible hyperk\"ahler manifolds, we can define an analogue of the positive cone of an irreducible hyperk\"ahler manifold. 

Recall that if $X$ is irreducible hyperk\"ahler and $\mathcal{C}_{X}$ is the cone of real $(1,1)-$classes over which the Beauville form is strictly positive, the positive cone $\mathcal{C}_{X}^{+}$ is the connected component of $\mathcal{C}_{X}$ which contains the K\"ahler cone of $X$. A result of Huybrechts shows that $\mathcal{C}_{X}^{+}$ is contained in (the interior of) the pseudo-effective cone of $X$. 

Theorem \ref{thm:thm2} is used to prove that on a compact, connected holomorphically symplectic $b_{2}-$manifold $X$ which is limit of irreducible hyperk\"ahler manifolds the intersection of the pseudo-effective cone of $X$ and of $\mathcal{C}_{X}$ (which can be defined as for irreducible hyperk\"ahler manifolds by Theorem \ref{thm:thm1}) consists of exactly one of the two connected components of $\mathcal{C}_{X}$: such a component is the positive cone of $X$, still denoted $\mathcal{C}_{X}^{+}$. We then prove:

\begin{thm}
\label{thm:thm3}Let $X$ be a compact, connected holomorphically symplectic $b_{2}-$manifold which is limit of irreducible hyperk\"ahler manifolds. If there is $\alpha\in\mathcal{C}_{X}^{+}$ such that 
\begin{enumerate}
 \item $\alpha\cdot C>0$ for every rational curve $C$ on $X$, and
 \item for every $\beta\in H^{1,1}(X)\cap H^{2}(X,\mathbb{Z})$ we have $q_{X}(\alpha,\beta)\neq 0$,
\end{enumerate}
then $X$ is irreducible hyperk\"ahler, and $\alpha$ is a K\"ahler class on $X$.
\end{thm}

The proof of this result is based on Theorem \ref{thm:thm2}, which gives a bimeromorphism $f:Y\dashrightarrow X$ between $X$ and an irreducible hyperk\"ahler manifold $Y$. Using twistor lines for real $(1,1)-$classes on $X$ (which can be defined similarily to the hyperk\"ahler case thanks to Theorem \ref{thm:thm1}) and a strategy used by Huybrechts for irreducible hyperk\"ahler manifolds, we show that the conditions on $\alpha$ imply that $f^{*}\alpha$ is a K\"ahler class on $Y$. An easy argument then shows that $f$ is a biholomorphism, and that $\alpha$ is a K\"ahler class on $X$.

The last part of the paper is devoted to show that on $M$ a class $\alpha$ as in the statement of Theorem \ref{thm:thm3} exists. This is obtained by using the (Hodge) isometry $$\lambda_{v}:v^{\perp}\otimes\mathbb{R}\longrightarrow H^{2}(M,\mathbb{R})$$(whose existence was proved in \cite{PT}) to produce classes in $\mathcal{C}_{M}$. By deforming to a moduli space of slope-stable sheaves on a projective K3 surface, and by using a classical construction of ample line bundles on $M$ in this case (starting from an ample line bundle on $S$), we are able to conclude that a class as in Theorem \ref{thm:thm3} exists, concluding the proof of Theorem \ref{thm:main}.

\subsection{The Beauville form and the Local Torelli Theorem}

The starting point of the proof of Theorem \ref{thm:main} is Theorem \ref{thm:thm1}, which is due to Guan. We will not prove it here (the proof can be found in \cite{Gua3}), but we recall the definition of $q_{X}$ and the Local Torelli Theorem.

\subsubsection{The Beauville form on $H^{2}(X,\mathbb{C})$.} Let $X$ be a compact, connected holomorphically symplectic manifold of complex dimension $2n$.

The \textit{Beauville form} of $X$ is a quadratic form on $H^{2}(X,\mathbb{C})$ defined as follows. First, choose a holomorphic symplectic form $\sigma$ on $X$, and assume for simplicity that $\int_{X}\sigma^{n}\wedge\overline{\sigma}^{n}=1$. For every $\alpha\in H^{2}(X,\mathbb{C})$, we let
$$q_{\sigma}(\alpha):=\frac{n}{2} \int_{X}\alpha^{2}\wedge\sigma^{n-1}\wedge\overline{\sigma}^{n-1}+(1-n)\int_{X}\alpha\wedge\sigma^{n}\wedge\overline{\sigma}^{n-1}\int_{X}\alpha\wedge\sigma^{n-1}\wedge\overline{\sigma}^{n}.$$Note that $q_\sigma(\sigma+\overline{\sigma})= (\int_X\sigma^n\wedge\overline{\sigma}^n)^2=1$ so $q_\sigma$ is non-trivial. Moreover, the quadratic form $q_{\sigma}$ depends a priori on the choice of $\sigma$.

\subsubsection{The period map.} Let $X$ be a compact, connected holomorphically symplectic $b_{2}-$manifold of complex dimension $2n$, and suppose that $h^{2,0}(X)=1$. We let $f:\mathcal{X}\longrightarrow B$ be its Kuranishi family, and $0\in B$ a point such that the fiber $X_{0}$ is isomorphic to $X$.

By Theorem 1 and the following Remark 1 of \cite{Gua3}, it follows that $B$ is smooth, and that up to shrinking it we can even suppose that all the fibers of the Kuranishi family are holomorphically symplectic.

Up to shrinking $B$, for every $b\in B$ the fiber $X_{b}$ of $f$ is a compact, connected holomorphically symplectic $b_{2}-$manifold (since being a $b_{2}-$manifold is an open property). Moreover, again up to shrinking $B$, by the Ehresmann Fibration Theorem we can suppose that $\mathcal{X}$ is diffeomorphic to $X\times B$. In particular, this induces a diffeomorphism $u_{b}:X\longrightarrow X_{b}$ for every $b\in B$, and hence an isomorphism of complex vector spaces $$u_{b}^{*}:H^{2}(X_{b},\mathbb{C})\longrightarrow H^{2}(X,\mathbb{C}).$$

We now let $\mathbb{P}:=\mathbb{P}(H^{2}(X,\mathbb{C}))$, and $$p:B\longrightarrow\mathbb{P},\,\,\,\,\,\,\,p(b):=[u_{b}^{*}(\sigma_{b})],$$where $\sigma_{b}$ is the holomorphic symplectic form on $X_{b}$ such that $\int_{X_{b}}\sigma_{b}^{n}\wedge\overline{\sigma}_{b}^{n}=1$ (we notice that such a $\sigma_{b}$ is unique as $h^{2,0}(X)=1$, and hence $h^{2,0}(X_{b})=1$). The map $p$ is holomorphic, and will be called \textit{period map of} $X$.

We let $Q_{\sigma}$ be the quadric defined by the quadratic form $q_{\sigma}$ in $\mathbb{P}$, i. e. $$Q_{\sigma}=\{\alpha\in\mathbb{P}\,|\,q_{\sigma}(\alpha)=0\},$$and $\Omega_{\sigma}$ be the open subset of $Q_{\sigma}$ defined as $$\Omega_{\sigma}:=\{\alpha\in Q_{\sigma}\,|\,q_{\sigma}(\alpha+\overline{\alpha})>0\}.$$ 

As showed in \cite{Gua3} we have the following, known as \textit{Local Torelli Theorem}:

\begin{prop}
\label{prop:tor}Let $X$ be a compact, connected holomorphically symplectic $b_{2}-$manifold, such that $h^{2,0}(X)=1$.
\begin{enumerate}
 \item The quadratic form $q_{\sigma}$ (and hence $Q_{\sigma}$ and $\Omega_{\sigma}$) is independent of $\sigma$, and will hence be denoted $q_{X}$ (and similarly $Q_{X}$ and $\Omega_{X}$).
 \item Up to a positive rational multiple, the quadratic form $q_{X}$ is a non-degenerate quadratic form on $H^{2}(X,\mathbb{Z})$ of signature $(3,b_{2}(X)-3)$.
 \item If $B$ is the base of the Kuranishi family of $X$, we have that $p(B)\subseteq\Omega_{X}$, and that $p:B\longrightarrow\Omega_{X}$ is a local biholomorphism.
\end{enumerate}
\end{prop}

Using the non-generate quadratic form $q_{X}$ (of signature $(3,b_{2}(X)-3)$) we let $$\mathcal{C}'_{X}:=\{\alpha\in H^{2}(X,\mathbb{R})\,|\,q_{X}(\alpha)>0\},$$which is an open cone in $H^{2}(X,\mathbb{R})$ having two connected components. Moreover, we let $$\widetilde{H}^{1,1}_{\mathbb{R}}(X):=\mathrm{Im}(\{\alpha\in H^{1,1}(X)\,|\,d\alpha=0\}\longrightarrow H^{2}(X,\mathbb{C}))\cap H^{2}(X,\mathbb{R}),$$and notice that this consists exactly of de Rham cohomology classes of real $d-$closed $(1,1)-$forms on $X$. We let $$\mathcal{C}_{X}:=\mathcal{C}'_{X}\cap\widetilde{H}^{1,1}_{\mathbb{R}}(X),$$which is then an open cone in $\widetilde{H}^{1,1}_{\mathbb{R}}(X)$ having two connected components. 

\section{Limits of irreducible hyperk\"ahler manifolds}

This section is devoted to prove that every compact, connected holomorphically symplectic $b_{2}-$manifold $X$ which is limit of irreducible hyperk\"ahler manifold, is bimeromorphic to an irreducible hyperk\"ahler manifold: in other words, we prove Theorem \ref{thm:thm2}.

The proof is divided in several sections. First we construct a moduli space $\mathcal{M}_{Z}$ of marked manifolds, and thanks to the Local Torelli Theorem we may give it the structure of a (non-separated) complex space. It will carry a period map to some period domain, which is locally a biholomorphism.

Then we show that each point in the closure of the open subset of $\mathcal{M}_{Z}$ given by irreducible hyperk\"ahler manifolds is non-separated from an irreducible hyperk\"ahler manifold.

Adapting an argument of Siu (for K3 surfaces) and Huybrechts (for higher dimensional irreducible hyperk\"ahler manifolds), we conclude the proof of Theorem \ref{thm:thm2}.

\subsection{The moduli space of $\Lambda-$marked manifolds}

In this section, we let $Z$ be an irreducible hyperk\"ahler manifold, and we write $(\Lambda,q):=(H^{2}(Z,\mathbb{Z}),q_{Z})$ for the Beauville lattice of $Z$. We let $\mathbb{P}_{\Lambda}:=\mathbb{P}(\Lambda\otimes\mathbb{C})$, and inside of it we let $$Q_{\Lambda}:=\{\alpha\in\mathbb{P}_{\Lambda}\,|\,q(\alpha)=0\},$$which is the quadric defined by $q$, and $$\Omega_{\Lambda}:=\{\alpha\in Q_{\Lambda}\,|\,q(\alpha+\overline{\alpha})>0\}.$$

The following is immediate, as $X$ is a $b_{2}-$manifold and the Hodge numbers are upper-semicontinuous.

\begin{prop}
\label{prop:h20}Let $X$ be a compact, connected holomorphically symplectic $b_{2}-$manifold which is deformation equivalent to an irreducible hyperk\"ahler manifold $Z$. Then for every $p,q\geq 0$ such that $p+q=2$ we have $h^{p,q}(X)=h^{p,q}(Z)$. In particular $h^{2,0}(X)=1$.
\end{prop}

If $X$ is a compact, connected holomorphically symplectic $b_{2}-$manifold which is deformation equivalent to $Z$, by Proposition \ref{prop:h20} and Theorem \ref{thm:thm1} we know that $H^{2}(X,\mathbb{Z})$ carries a non-degenerate quadratic form $q_{X}$, and that there is an isometry $\phi:H^{2}(X,\mathbb{Z})\longrightarrow\Lambda$. The isometry $\phi$ is called $\Lambda-$\textit{marking} on $X$, and the pair $(X,\phi)$ is a $\Lambda-$\textit{marked manifold}. The set of $\Lambda-$marked manifolds will be denoted $\mathcal{M}'_{Z}$.

Moreover, we let $\mathcal{M}_{Z}:=\mathcal{M}'_{Z}/\sim$, where $(X,\phi)\sim(X',\phi')$ if and only if there is a biholomorphism $f:X\longrightarrow X'$ such that $f^{*}\circ\phi'=\phi$. The set $\mathcal{M}_{Z}$ will be referred to as the \textit{moduli space of} $\Lambda-$\textit{marked manifolds}. We let $\mathcal{M}^{hk}_{Z}$ be the subset of $\mathcal{M}_{Z}$ of pairs $(X,\phi)$ where $X$ is an irreducible hyperk\"ahler manifold: it will be called \textit{moduli space of} $\Lambda-$\textit{marked hyperk\"ahler manifolds}. 

We first show that $\mathcal{M}_{Z}$ has the structure of complex space (hence justifying the name \textit{space} we use for it): the following is a generalization of Proposition 4.3 of \cite{H2}, and requires the same proof.

\begin{prop}
\label{prop:mlambda}Let $Z$ be an irreducible hyperk\"ahler manifold and $(\Lambda,q)$ its Beauville lattice.
\begin{enumerate}
 \item For any $(X,\phi)\in\mathcal{M}_{Z}$ there is an inclusion $i_{X}:B\longrightarrow\mathcal{M}_{Z}$, where $B$ is the base of the Kuranishi family of $X$.
 \item The set $\mathcal{M}_{Z}$ has the structure of smooth complex space of dimension $b_{2}(Z)-2$. 
 \item The subset $\mathcal{M}^{hk}_{Z}$ is an open subset (in the analytic topology) of $\mathcal{M}_{Z}$.
\end{enumerate}
\end{prop}

\proof Let $X$ be a compact, connected holomorphically symplectic $b_{2}-$manifold which is deformation equivalent to $Z$, and $f:\mathcal{X}\longrightarrow B$ its Kuranishi family. 

Up to shrinking $B$ we can suppose that it is a complex disk of dimension $b_{2}(X)-2=b_{2}(Z)-2$, and as we have seen before for every $b\in B$ we can suppose that $X_{b}$ is a compact, connected holomorphically symplectic $b_{2}-$manifold (which is clearly deformation equivalent to $Z$). 

Moreover, we can suppose that $\mathcal{X}$ is diffeomorphic (over $B$) to the trivial family $X\times B$, and that we have a diffeomorphism $u_{b}:X\longrightarrow X_{b}$ inducing an isometry $u_{b}^{*}:H^{2}(X_{b},\mathbb{Z})\longrightarrow H^{2}(X,\mathbb{Z})$. We let $\phi_{b}:=\phi\circ u_{b}^{*}$, which is then a $\Lambda-$marking on $X_{b}$, for every $b\in B$. 

It follows that for every $b\in B$ we have $(X_{b},\phi_{b})\in\mathcal{M}_{Z}$, so that we have a map $$i_{X}:B\longrightarrow\mathcal{M}_{Z},\,\,\,\,\,\,\,\,\,\,i_{X}(b):=(X_{b},\phi_{b}).$$

We show that $i_{X}$ is an inclusion. Let $b,b'\in B$ and suppose that $i_{X}(b)=i_{X}(b')$. This means that $(X_{b},\phi_{b})\sim(X_{b'},\phi_{b'})$, i. e. there is a biholomorphism $f:X_{b}\longrightarrow X_{b'}$ such that $$f^{*}=\phi_{b}^{-1}\circ\phi_{b'}.$$By definition of $\phi_{b}$ and $\phi_{b'}$, this means that $$f^{*}=(\phi\circ u_{b}^{*})^{-1}\circ(\phi\circ u_{b'}^{*})=(u_{b}^{*})^{-1}\circ u_{b'}^{*}.$$

Now, let $\sigma_{b}$ and $\sigma_{b'}$ be symplectic forms on $X_{b}$ and $X_{b'}$ respectively. As $f$ is a biholomorphism, the form $f^{*}\sigma_{b'}$ is holomorphic symplectic on $X_{b}$, and hence $$[u_{b}^{*}\sigma_{b}]=[u_{b}^{*}f^{*}\sigma_{b'}].$$But as $f^{*}=(u_{b}^{*})^{-1}\circ u_{b'}^{*}$, this implies that $[u_{b}^{*}\sigma_{b}]=[u_{b'}^{*}\sigma_{b'}]$. By definition of the period map of $X$, this means that $p(b)=p(b')$.

But now recall that by point (3) of Proposition \ref{prop:tor}, the period map $p:B\longrightarrow\Omega$ is a local biholomorphism: up to shrinking $B$, for $b\neq b'\in B$ we have $p(b)\neq p(b')$. It follows that up to shrinking $B$ the condition $i_{X}(b)=i_{X}(b')$ implies $b=b'$, and $i_{X}$ is an inclusion of $B$ in $\mathcal{M}_{Z}$. This proves point 1 of the statement.

To give $\mathcal{M}_{Z}$ the structure of a complex space, we just need to show that each point of $\mathcal{M}_{Z}$ has a neighborhood having the structure of a complex manifold, and that whenever two neighborhoods of this type intersect, the corresponding complex structures glue.

If $(X,\phi)\in\mathcal{M}_{Z}$, the previous part of the proof suggests to view $i_{X}(B)$ as a neighborhood $(X,\phi)$ in $\mathcal{M}_{Z}$. Now, let $(X,\phi),(X',\phi')\in\mathcal{M}_{Z}$, and we let $B$ and $B'$ be the bases of the Kuranishi families of $X$ and of $X'$, respectively. If $i_{X}(B)\cap i_{X}(B')\neq\emptyset$, then $B\cap B'$ is an open subset of $B$ and $B'$, over which the Kuranishi families coincide. This allows us to glue the Kuranishi families along $B\cap B'$, and hence the complex structures of $i_{X}(B)$ and $i_{X}(B')$ can be glued in $\mathcal{M}_{Z}$. This shows that $\mathcal{M}_{Z}$ has the structure of a complex space.

We notice that as each base $B$ of a Kuranishi family of a compact, connected holomorphically symplectic $b_{2}-$manifold is smooth (see section 2.2) of dimension $b_{2}(Z)-2$, it follows that $\mathcal{M}_{Z}$ is a smooth complex space, and its dimension is $b_{2}(Z)-2$. This proves point 2 of the statement.

The fact that $\mathcal{M}^{hk}_{\Lambda}$ is open in the analytic topology is a consequence of the fact that K\"ahlerness is an open property in the analytic topology.\endproof

The complex space $\mathcal{M}_{Z}$ has two connected components, and one can pass from one to the other by sending $(X,\phi)$ to $(X,-\phi)$.

We now define the period map in this generality: we let $$\pi:\mathcal{M}_{Z}\longrightarrow\mathbb{P}_{\Lambda},\,\,\,\,\,\,\,\pi(X,\phi):=[\phi(\sigma)],$$where $\sigma$ is a holomorphic symplectic form on $X$. Notice that the $\Lambda-$marking $\phi$ induces an isomorphism $$\overline{\phi}:\mathbb{P}\longrightarrow\mathbb{P}_{\Lambda},$$and as it is an isometry it induces an isomorphism $$\overline{\phi}:\Omega_{X}\longrightarrow\Omega_{\Lambda}.$$If $B$ is the base of the Kuranishi family of $X$, we have $\pi_{|i_{X}(B)}=\overline{\phi}\circ p$: if $b\in B$ and $\sigma_{b}$ is a symplectic form on $X_{b}$, we have $$\overline{\phi}(p(b))=\overline{\phi}[u_{b}^{*}\sigma_{b}]=[\phi(u_{b}^{*}\sigma_{b})]=[\phi_{b}(\sigma_{b})]=\pi(b).$$

The first two points of the following Proposition are just a translation in this language of Theorem \ref{thm:thm1}. For the last point: the surjectivity is Theorem 8.1 of \cite{Hu99}; the general injectivity is the Global Torelli Theorem of Verbistky.

\begin{prop}
\label{prop:periods}We have the following properties:
\begin{enumerate}
 \item the image of $\pi$ is contained in $\Omega_{\Lambda}$;
 \item the map $\pi$ is a local biholomorphism;
 \item if $\mathcal{M}^{hk,0}_{Z}$ is a connected component of $\mathcal{M}^{hk}_{Z}$, the map $\pi_{|\mathcal{M}^{hk,0}_{Z}}$ is surjective and generically injective.
\end{enumerate}
\end{prop}

Now, we let $\overline{\mathcal{M}}^{hk}_{Z}$ be the closure of $\mathcal{M}^{hk}_{Z}$ in $\mathcal{M}_{Z}$. Using this formalism, we can state Theorem \ref{thm:thm2} in an equivalent way:

\begin{prop}
\label{prop:thm2equiv}
If $(X,\phi)\in\overline{\mathcal{M}}^{hk}_{Z}$, then $X$ is bimeromorphic to an irreducible hyperk\"ahler manifold (hence, it is in the Fujiki class $\mathcal{C}$).
\end{prop}

This is the statement we will prove in the next sections.

\subsection{Non-separatedness in $\mathcal{M}_{Z}$}

The first result we show is the following:

\begin{prop}
\label{prop:nonsep}Let $(X,\phi)\in\overline{\mathcal{M}}^{hk,0}_{Z}$. Then there is $(Y,\psi)\in\mathcal{M}^{hk,0}_{Z}$ such that $(X,\phi)$ and $(Y,\psi)$ are non-separated in $\mathcal{M}_{Z}$.
\end{prop}

\proof The statement is trivial if $(X,\phi)\in\mathcal{M}^{hk,0}_{Z}$. We then suppose that $(X,\phi)\in\overline{\mathcal{M}}^{hk,0}_{Z}\setminus\mathcal{M}^{hk,0}_{Z}$.

We let $p:=\pi(X,\phi)\in\Omega_{\Lambda}$ be the period of $(X,\phi)$. As $\pi_{|\mathcal{M}^{hk,0}_{Z}}$ is surjective, there is $(Y,\psi)\in\mathcal{M}^{hk,0}_{Z}$ such that $\pi(Y,\psi)=p$. We show that $(X,\phi)$ and $(Y,\psi)$ are non-separated in $\mathcal{M}_{Z}$.

To do so, let $U_{X}$ and $U_{Y}$ be two open neighborhoods of $(X,\phi)$ and $(Y,\psi)$ respectively in $\mathcal{M}_{Z}$. Up to shrinking $U_{X}$ and $U_{Y}$, we can suppose that $\pi(U_{X})=\pi(U_{Y})=:V$. 

Moreover, by point (3) of Proposition \ref{prop:tor}, up to shrinking $U_{X}$ and $U_{Y}$ we can suppose that $\pi_{|U_{Y}}:U_{Y}\longrightarrow V$ and $\pi_{|U_{X}}:U_{X}\longrightarrow V$ are biholomorphisms. Finally, as K\"ahlerness is an open property in the analytic topology, up to shrinking $U_{X}$ and $U_{Y}$ we can suppose that $U_{Y}\subseteq\mathcal{M}^{hk,0}_{Z}$.

Now, as $(X,\phi)\in\overline{\mathcal{M}}^{hk,0}_{Z}$, there is a hyperk\"ahler manifold $X'$ and a marking $\phi'$ on $X'$ such that $(X',\phi')\in U_{X}\cap\mathcal{M}^{hk,0}_{Z}$. We can choose $(X',\phi')$ to be generic. Let $p':=\pi(X',\phi')\in V$: as $\pi_{|U_{Y}}:U_{Y}\longrightarrow V$ is surjective, there is $(Y',\psi')\in U_{Y}$ such that $\pi(Y',\psi')=p'$, and as $U_{Y}\subseteq\mathcal{M}^{hk,0}_{Z}$, we have that $Y'$ is an irreducible hyperk\"ahler manifold.

Hence $(X',\phi')$ and $(Y',\psi')$ are two generic points in $\mathcal{M}_{Z}^{hk,0}$: by point (3) of Proposition \ref{prop:periods} we then have $(X',\phi')=(Y',\psi')$ in $\mathcal{M}_{Z}$, so that $U_{X}\cap U_{Y}\neq\emptyset$, and we are done.\endproof

This result will be the starting point of the proof of Proposition \ref{prop:thm2equiv}.

\subsection{The proof of Theorem \ref{thm:thm2}}

We now prove a key result in the proof of Proposition \ref{prop:thm2equiv}.

\begin{lem}
\label{lem:border}Let $B$ be a connected complex manifold and $\mathcal{X}\longrightarrow B$ and $\mathcal{Y}\longrightarrow B$ be two smooth, proper families verifying the following properties:
\begin{enumerate}
 \item for every $b\in B$ the fiber $Y_{b}$ is an irreducible hyperk\"ahler manifold with a $\Lambda-$marking $\psi_{b}$;
 \item for every $b\in B$ the fiber $X_{b}$ is a compact, connected holomorphically symplectic $b_{2}-$manifold deformation equivalent to an irreducible hyperk\"ahler manifolds, which has a $\Lambda-$marking $\phi_{b}$; 
 \item there is a non-empty open subset $V$ of $B$ such that for each $b\in V$ there is an isomorphism $f_{b}:Y_{b}\longrightarrow X_{b}$ such that $f_{b}^{*}=\psi_{b}^{-1}\circ\phi_{b}$;
 \item for generic $b\in V$ we have $H^{1,1}(X_{b})\cap H^{2}(X_{b},\mathbb{Z})=0$.
\end{enumerate}
Then $V$ is dense in $B$.
\end{lem}

\proof We choose a point $0\in B$ and let $X:=X_{0}$ and $Y:=Y_{0}$. We show that $\partial V:=\overline{V}\setminus V$ is contained in a countable union of analytic subvarieties of $B$. It has then real codimension at least 2 in $B$, hence it cannot separate the disjoint open subsets $V$ and $B\setminus\overline{V}$. As $V\neq\emptyset$, it follows that $B=\overline{V}$.

In order to show that $\partial V$ is contained in a countable union of analytic subvarieties of $B$, we show that if $s\in\partial V$, then $Y_{s}$ has either effective divisors or curves. This implies that $$s\in\bigcup_{\alpha}S_{\alpha},$$where $\alpha\in H^{2}(Y,\mathbb{Z})$ and $S_{\alpha}$ is the analytic subvariety of $B$ given by those $b\in B$ such that $\alpha\in NS(Y_{b})$. We proceed by contradiction: we let $s\in\partial V$, and we suppose that $Y_{s}$ has no effective divisors and no curves. 

As $s\in\partial V$, it follows that $(X_{s},\phi_{s})$ and $(Y_{s},\psi_{s})$ are non-separated points in $\mathcal{M}_{Z}$, were $Z$ is an irreducible hyperk\"ahler manifold among all the $Y_{b}$. We first show that $X_{s}$ and $Y_{s}$ are bimeromorphic. To show this, let $\beta_{s}$ be a K\"ahler form on $Y_{s}$ and $\alpha_{s}$ a closed real $(1,1)-$form on $X_{s}$ whose cohomology class is in $\mathcal{C}_{X_{s}}$. 

Consider a continuous family $\{\beta_{t}\}_{t\in B}$, where $\beta_{t}$ is a closed $(1,1)-$form on $Y_{t}$. As K\"ahlerness is an open property in the analytic topology, there is an analytic open neighborhood $U$ of $s$ in $B$ such that for every $t\in U$ the form $\beta_{t}$ is K\"ahler on $Y_{t}$.

Moreover, consider a continuous family $\{\alpha_{t}\}_{t\in B}$ where $\alpha_{t}$ is a closed $(1,1)-$form on $X_{t}$. Up to shrinking $U$, and as the positivity of $q_{X}$ is an open property, we can suppose that for every $t\in U$ the cohomology class of $\alpha_{t}$ is in $\mathcal{C}_{X_{t}}$.

Notice that as $s\in\partial V$, the intersection of $U$ with $V$ is not empty, and the generic point $t\in U\cap V$ is such that $f_{t}:Y_{t}\longrightarrow X_{t}$ is a biholomorphism such that $f_{t}^{*}=\psi_{t}^{-1}\circ\phi_{t}$. By hypothesis, we have that $X_{t}$ is an irreducible hyperk\"ahler manifold with $NS(X_{t})=0$. By Corollary 5.7 of \cite{Hu99}, this implies that the K\"ahler cone of $X_{t}$ is one of the two connected components of $\mathcal{C}_{X_{t}}$.

As the cohomology class $[\alpha_{t}]$ of $\alpha_{t}$ is in $\mathcal{C}_{X_{t}}$, it follows that either $[\alpha_{t}]$ or $-[\alpha_{t}]$ is in the K\"ahler cone of $X_{t}$. Up to changing the sign of $\alpha_{t}$, we can then suppose that for the generic $t\in U\cap V$ the class $[\alpha_{t}]$ is a K\"ahler class, and that $\alpha_{t}$ is a K\"ahler form.

In conclusion, there is a sequence $\{t_{m}\}_{m\in\mathbb{N}}$ of point of $U\cap V$ which converges to $s$, and such that for every $m\in\mathbb{N}$ we have a K\"ahler form $\alpha_{m}:=\alpha_{t_{m}}$ on $X_{m}:=X_{t_{m}}$ and a K\"ahler form $\beta_{m}:=\beta_{t_{m}}$ on $Y_{m}:=Y_{t_{m}}$, such that $\alpha_{m}$ converges to $\alpha_{s}$ and $\beta_{m}$ converges to $\beta_{s}$.

As $t_{m}\in V$, we have an isomorphism $f_{m}:Y_{m}\longrightarrow X_{m}$ such that $f_{m}^{*}=\psi_{m}^{-1}\circ\phi_{m}$ (where $\psi_{m}:=\psi_{t_{m}}$ and $\phi_{m}=\phi_{t_{m}}$). We let $\Gamma_{m}$ be the graph of $f_{m}$, and we compute its volume in $X_{m}\times Y_{m}$ with respect to the the form $p_{1}^{*}\alpha_{m}+p_{2}^{*}\beta_{m}$, where $p_{1}$ and $p_{2}$ are the projections of $X_{m}\times Y_{m}$ to $X_{m}$ and $Y_{m}$ respectively.

We then have $$\vol(\Gamma_{m})=\int_{Y_{m}}(\beta_{m}+f_{m}^{*}\alpha_{m})^{2n}=\int_{Y_{m}}([\beta_{m}]+f_{m}^{*}[\alpha_{m}])^{2n}=$$ $$=\int_{Y_{m}}([\beta_{m}]+\psi_{m}^{-1}\circ\phi_{m}([\alpha_{m}]))^{2n}.$$Taking the limit for $m$ going to infinity we get $$\lim_{m\rightarrow+\infty}\vol(\Gamma_{m})=\int_{Y_{s}}([\beta_{s}]+\psi_{s}^{-1}\circ\phi_{s}([\alpha_{s}]))^{2n}<+\infty.$$

Hence, the volumes of the $\Gamma_{m}$ are bounded, so that by the Bishop Theorem the cycles $\Gamma_{m}$ converge to a cycle $\Gamma$ of $X_{s}\times Y_{s}$ with the same cohomological properties of the $\Gamma_{m}$'s: namely, we have $[\Gamma]\in H^{4n}(X_{s}\times Y_{s},\mathbb{Z})$, and if $p_{1}$ and $p_{2}$ are the two projections from $X_{s}\times Y_{s}$ to $X_{s}$ and $Y_{s}$ respectively, we have $p_{1*}[\Gamma]=[X_{s}]$, $p_{2*}[\Gamma]=[Y_{s}]$, and $$[\Gamma]_{*}\gamma:=p_{2*}([\Gamma]\cdot p_{1}^{*}\gamma)=\psi_{s}^{-1}(\phi_{s}(\gamma)),$$for every $\gamma\in H^{2}(X_{s},\mathbb{Z})$.

Now, let us split $\Gamma$ in its irreducible components: by the previous properties we then have that either $\Gamma=Z+\sum_{i}D_{i}$ where $p_{1}:Z\longrightarrow X_{s}$ and $p_{2}:Z\longrightarrow Y_{s}$ are both generically one-to-one, or $\Gamma=Z_{1}+Z_{2}+\sum_{i}D_{i}$ where $p_{1}:Z_{1}\longrightarrow X_{s}$ and $p_{2}:Z_{2}\longrightarrow Y_{s}$ are generically one-to-one, but neither $p_{1}:Z_{2}\longrightarrow Y_{s}$ nor $p_{2}:Z_{1}\longrightarrow X_{s}$ are generically finite. In both cases we have $p_{1*}[D_{i}]=p_{2*}[D_{i}]=0$. As shown in \cite{H2} (proof of Theorem 4.3), the second case cannot happen: it follows that $Z$ is a bimeromorphism between $X_{s}$ and $Y_{s}$.

Now, recall that $Y_{s}$ is supposed to have no effective divisors nor non-trivial curves. It follows from this that the bimeromorphism $Z$ is indeed an isomorphism: hence there is an isomorphism $f:Y_{s}\longrightarrow X_{s}$ whose graph is $Z$, and $X_{s}$ is K\"ahler (as $Y_{s}$ is). Moreover, by Corollary 5.7 of \cite{Hu99}, the K\"ahler cone of $Y_{s}$, and hence that of $X_{s}$, is one of the components of $\mathcal{C}_{Y_{s}}$, and the morphisms $[D_{i}]_{*}:H^{2}(X_{s},\mathbb{Z})\longrightarrow H^{2}(Y_{s},\mathbb{Z})$ are all trivial (see Lemma 5.5 in \cite{Hu99}). In particular, we get that $f^{*}=[Z]_{*}=[\Gamma]_{*}$

We now prove that $\Gamma=Z$. To do this, recall that the class of $\alpha_{s}$ is K\"ahler, and let $\gamma_{s}:=f^{*}\alpha_{s}$, which is K\"ahler again on $Y_{s}$. We compute the volumes on $X_{s}\times Y_{s}$ with respect to the K\"ahler class $p_{1}^{*}\alpha_{s}+p_{2}^{*}\gamma_{s}$. We have $$\vol(\Gamma)=\vol(Z)+\sum_{i}\vol(D_{i})=\int_{Z}(p_{1}^{*}\alpha_{s}+p_{2}^{*}\gamma_{s})^{2n}+\sum_{i}\vol(D_{i})=$$ $$=\int_{Y_{s}}([Z]_{*}\alpha_{s}+\gamma_{s})^{2n}+\sum_{i}\vol(D_{i})=\int_{Y_{s}}(f^{*}\alpha_{s}+\gamma_{s})^{2n}+\sum_{i}\vol(D_{i})=$$ $$=\int_{Y_{s}}(2\gamma_{s})^{2n}+\sum_{i}\vol(D_{i}),$$where we used that $f_{*}=[Z]_{*}=[\Gamma]_{*}$.

Now, choose a sequence $\{t_{m}\}_{m\in\mathbb{N}}$ of points in $V\cap U$ converging to $s$, and consider the graphs $\Gamma_{m}$ of the isomorphisms $f_{m}:Y_{m}\longrightarrow Y_{m}$. Letting $\gamma_{m}:=f_{m}^{*}\alpha_{m}$, the sequence $\{\gamma_{m}\}$ converges to $\gamma_{s}$, and the sequence $\{\vol(\Gamma_{m})\}_{m\in\mathbb{N}}$ converges to $\vol(\Gamma)$. But $$\vol(\Gamma_{m})=\int_{Y_{m}}([\Gamma_{m}]_{*}\alpha_{m}+\gamma_{m})^{2n}=\int_{Y_{m}}(2\gamma_{m})^{2n},$$hence this sequence converges to $$\int_{Y_{s}}(2\gamma_{s})^{2n},$$so that $$\vol(\Gamma)=\int_{Y_{s}}(2\gamma_{s})^{2n}.$$It turns then out that $\vol(D_{i})=0$, hence $D_{i}=0$ and $\Gamma=Z$.

It follows that $$f^{*}=[Z]_{*}=[\Gamma]_{*}=\psi_{s}^{-1}\circ\phi_{s}.$$But this implies that $s\in V$, contradicting that $s\in\partial V$. In conclusion if $s\in\partial V$ on $Y_{s}$ there are either effective divisors or curves, and we are done.\endproof

We are now able to conclude the proof of Proposition \ref{prop:thm2equiv}.

\proof By Proposition \ref{prop:nonsep} we know that there is an irreducible symplectic manifold $Y$, together with a marking $\psi$, such that $(X,\phi)$ is non-separated from $(Y,\psi)$ in $\mathcal{M}_{Z}$. 

Consider $Def(X)$ and $Def(Y)$, the bases of the Kuranishi families of $X$ and $Y$ respectively. Up to shrinking them, as the points $(X,\phi)$ and $(Y,\psi)$ are non-separated, the Local Torelli Theorem allows us to identify them. Hence the Kuranishi families of $X$ and $Y$ are over the same base $B$, and we suppose that $X$ and $Y$ are over the same point $0\in B$.

The non-separatedness implies that there is $t\in B$ such that $X_{t}$ and $Y_{t}$ are isomorphic under an isomorphism $f_{t}$ such that $f_{t}^{*}=\psi_{t}^{-1}\circ\phi_{t}$. Let $V$ be the biggest open subset of $B$ given by all $t\in B$ verifying this same property. As $V$ is open in $B$, and as $B$ is open in $\mathcal{M}_{Z}$, the generic point $t$ of $V$ is such that $NS(X_{t})=0$. 

We can then apply Lemma \ref{lem:border} to conclude that $\overline{V}=B$. Now, notice that if $0\in V$, then $X$ and $Y$ are isomorphic, and we are done. Hence we can suppose that $0\notin V$, so that $0\in\partial V$. We can then apply the same argument in the proof of Lemma \ref{lem:border} to conclude that $X$ and $Y$ are bimeromorphic.\endproof

\section{Criterion for K\"ahlerness}

We now want to prove a K\"ahlerness criterion for a compact, connected holomorphic symplectic $b_{2}-$manifold $X$ which is limit of irreducible hyperk\"ahler manifolds. 

Let us first recall a notation. In this section, $X$ is a compact, connected holomorphic symplectic $b_{2}-$manifold $X$ which is limit of irreducible hyperk\"ahler manifolds. By Theorem \ref{thm:thm2}, we know that $X$ is in the Fujiki class $\mathcal{C}$, hence $H^{2}(X,\mathbb{C})$ has a Hodge decomposition. In particular, we have $$\widetilde{H}^{1,1}_{\mathbb{R}}(X)=H^{2}(X,\mathbb{R})\cap H^{1,1}(X)=:H^{1,1}(X,\mathbb{R}).$$

\subsection{Twistor lines}

If $\sigma$ is a holomorphic symplectic form on $X$, the cohomology class of $\sigma$ allows us to define a real plane $$P(X):=(\mathbb{C}\cdot\sigma\oplus\mathbb{C}\cdot\overline{\sigma})\cap H^{2}(X,\mathbb{R})$$in $H^{2}(X,\mathbb{R})$, which is independent of $\sigma$ (as $h^{2,0}(X)=1$). If $\alpha\in H^{1,1}(X,\mathbb{R})$ we let $$F(\alpha):=P(X)\oplus\mathbb{R}\cdot\alpha,$$which is a $3-$dimensional real subspace of $H^{2}(X,\mathbb{R})$, and we let $F(\alpha)_{\mathbb{C}}:=F(\alpha)\otimes\mathbb{C}$.

Now, let $Z$ be an irreducible hyperk\"ahler manifold which is deformation equivalent to $X$, and let $(\Lambda,q)$ its Beauville lattice. If $\phi:H^{2}(X,\mathbb{Z})\longrightarrow\Lambda$ is a $\Lambda-$marking on $X$ (which exists by Theorem \ref{thm:thm1}), consider the point $(X,\phi)$ of $\mathcal{M}_{Z}$. As $X$ is limit of irreducible hyperk\"ahler manifolds, we have $(X,\phi)\in\overline{\mathcal{M}}^{hk}_{Z}$. 

Notice that $F(\alpha)_{\mathbb{C}}$ is a $3-$dimensional linear subspace of $H^{2}(X,\mathbb{C})$, hence $\phi(F(\alpha)_{\mathbb{C}})$ is a $3-$dimensional subspace of $\Lambda\otimes\mathbb{C}$, and $\mathbb{P}(\phi(F(\alpha)_{\mathbb{C}})$ is a plane in $\mathbb{P}_{\Lambda}$. Hence $\mathbb{P}(\phi(F(\alpha)_{\mathbb{C}}))\cap\Omega_{\Lambda}$ is a curve in $\Omega_{\Lambda}$ passing through $\pi(X,\phi)$. 

If $B$ is the base of the Kuranishi family of $X$, the inverse image $$T(\alpha):=\pi^{-1}(\mathbb{P}(\phi(F(\alpha)_{\mathbb{C}}))\cap\Omega_{\Lambda})\cap B$$is then a curve in $B$, which will be called \textit{twistor line of} $\alpha$. The restriction of the Kuranishi family of $X$ to $T(\alpha)$ will be denoted $$\kappa_{\alpha}:\mathcal{X}(\alpha)\longrightarrow T(\alpha).$$For every $t\in T(\alpha)$ there is real $(1,1)-$class $\alpha_{t}$ on the fiber $X_{t}$ of the Kuranishi family of $X$ over $t$, and the sequence $\{\alpha_{t}\}$ converges to $\alpha$. If $\alpha$ is K\"ahler, then $T(\alpha)\simeq\mathbb{P}^{1}$, and all $\alpha_{t}$ are K\"ahler on $X_{t}$.

\subsection{Cones in $H^{1,1}(X,\mathbb{R})$}

Using the notation introduced in the previous section, we have $$\mathcal{C}_{X}=\{\alpha\in H^{1,1}(X,\mathbb{R})\,|\,q_{X}(\alpha)>0\},$$which is an open cone in $H^{1,1}(X,\mathbb{R})$ having two connected components. 

If $X$ is K\"ahler (and hence irreducible hyperk\"ahler), the \textit{K\"ahler cone} $\mathcal{K}_{X}$ of $X$ (i. e. the open convex cone of K\"ahler classes on $X$) is contained in one of them: such a component is usually called \textit{positive cone} of $X$, and denoted $\mathcal{C}^{+}_{X}$. The other component will be denoted $\mathcal{C}^{-}_{X}$. If $NS(X)=0$, Corollary 5.7 of \cite{Hu99} gives us that $\mathcal{K}_{X}=\mathcal{C}^{+}_{X}$, a fact that has already been used in the previous sections.

Theorem 1.1 of \cite{Bou} tells us that if $X$ is irreducible hyperk\"ahler, then $\alpha\in\mathcal{C}^{+}_{X}$ is in the K\"ahler cone of $X$ if and only if $$\int_{C}\alpha>0$$for every rational curve $C$ of $X$. Our aim is to show a similar result for a compact, connected holomorphic symplectic $b_{2}-$manifold $X$ which is limit of irreducible hyperk\"ahler manifolds.

As on such a manifold the K\"ahler cone could be empty, we cannot use it to define the positive cone of $X$. Anyway, we can use the \textit{pseudo-effective cone} $\mathcal{E}_{X}$ of $X$, i. e. the closed connected cone of classes of positive closed real $(1,1)-$currents on $X$. If $X$ is irreducible hyperk\"ahler, by point i) of Theorem 4.3 of \cite{B2} we have $\mathcal{C}_{X}^{+}\subseteq\mathcal{E}_{X}$.

Popovici and Ugarte (see Theorem 5.9 of \cite{PU}) showed that if $\mathcal{X}\longrightarrow B$ is a smooth and proper family of sGG manifolds and $\{b_{n}\}$ is a sequence of points of $B$ converging to a point $b\in B$, then the limit of the pseudo-effective cones of $X_{b_{n}}$ is contained in $\mathcal{E}_{X_{b}}$, i. e. the pseudo-effective cone varies upper-semicontinuously along $B$.

As all manifolds in the Fujiki class $\mathcal{C}$ are sGG manifolds (see \cite{PU}), we conclude that the pseudo-effective cone varies upper-semicontinuously in families of class $\mathcal{C}$ manifolds. 

We now prove the following general fact about convex cones in a real finitely dimensional vector space.

\begin{lem}
\label{lem:coni}Let $V$ a real vector space of finite dimension $n$, and let $A,B\subseteq V$ two cones in $V$ such that:
\begin{enumerate}
 \item the cone $A$ is strictly convex (i. e. it does not contain any linear subspace of $V$) and closed;
 \item the cone $B$ is open and has two connected components, each of which is convex;
 \item for every $a\in B$, we have either $a\in A$ or $-a\in A$.
\end{enumerate}
Then $A\cap B$ is one of the connected components of $B$.
\end{lem}

\proof We first notice that if $B^{+}$ and $B^{-}$ are the two connected components of $B$, if $B^{+}\subseteq A$ we have $B^{+}=B\cap A$. Indeed, if $b'\in B^{-}\cap A$, then $-b'\in B^{+}\subseteq A$. It follows that $b',-b'\in A$, which is not possible as $A$ is a strictly convex cone.

We are left to prove that there is a connected component of $B$ which is contained in $A$. To do so, let $b_{0}\in B\cap A$, and let $B^{+}$ be the connected component of $B$ which contains $b_{0}$. We show that if $b_{1}\in B^{+}$, then $b_{1}\in A$.

Consider the segment $$[b_{0},b_{1}]:=\{b_{t}:=(1-t)b_{0}+tb_{1}\,|\,t\in[0,1]\}.$$Suppose that $b_{1}\notin A$, we have to find a contradiction. First, notice that as $b_{1}\notin A$, there is $t\in[0,1)$ such that $b_{t}\notin A$: indeed, if for every $t\in[0,1)$ we have $b_{t}\in A$, as $A$ is closed we would have $b_{1}\in A$.

As $b_{t}$ and $b_{1}$ are not in $A$, for every $s\in[t,1]$ we have $b_{s}\notin A$: indeed, as $b_{t},b_{1}\notin A$ but $b_{t},b_{1}\in B$, we have $-b_{t},-b_{1}\in A$. As $A$ is convex, the segment $[-b_{t},-b_{1}]$ (whose elements are the $-b_{s}$ for $s\in[t,1]$) is contained in $A$. But this means that $b_{s}\notin A$ as $A$ is a strictly convex cone.

The set of those $t\in[0,1]$ for which $b_{t}\notin A$ surely has an infimum $t_{0}\in[0,1]$. Hence, for every $t<t_{0}$ we have $b_{t}\in A$, and for every $t>t_{0}$ we have $b_{t}\notin A$. As $b_{t}\in B$, this implies that $-b_{t}\in A$ for every $t>t_{0}$. But as $A$ is closed, these conditions give $b_{t_{0}}\in A$ (as $b_{t}\in A$ for every $t<t_{0}$) and $-b_{t_{0}}\in A$ (as $-b_{t}\in A$ for every $t>t_{0}$). But as $A$ is a strictly convex cone, we get a contradiction.\endproof

This fact will be used in the proof of the following:

\begin{lem}
\label{lem:poseff}Let $X$ be a compact, connected holomorphic symplectic $b_{2}-$manifold which is limit of irreducible hyperk\"ahler manifolds. Then $\mathcal{C}_{X}\cap\mathcal{E}_{X}$ consists of exactly one connected component of $\mathcal{C}_{X}$.
\end{lem}

\proof The pseudo-effective cone $\mathcal{E}_{X}$ is strictly convex and closed in $H^{1,1}(X,\mathbb{R})$. The cone $\mathcal{C}_{X}$ is open ans has two connected components, each of which is convex. We show that if $\alpha\in\mathcal{C}_{X}$, then either $\alpha\in\mathcal{E}_{X}$ or $-\alpha\in\mathcal{E}_{X}$ (which in particular implies that $\mathcal{C}_{X}\cap\mathcal{E}_{X}\neq\emptyset$). Once this is done, the statement follows from Lemma \ref{lem:coni}.

Fix an irreducible hyperk\"ahler manifold $Z$ which is deformation equivalent to $X$, and we let $(\Lambda,q)$ be its Beauville lattice. Moreover, let $\alpha\in\mathcal{C}_{X}$, and consider a $\Lambda-$marking $\phi$ on $X$ (whose existence comes from Theorem \ref{thm:thm1}). As $X$ is limit of irreducible hyperk\"ahler manifolds, we have $(X,\phi)\in\overline{\mathcal{M}}^{hk}_{Z}$.

Let $\mathcal{X}\longrightarrow B$ be the Kuranishi family of $X$, and let $0$ be the point of $B$ over which the fiber of $X_{0}$ is $X$. By Theorem \ref{thm:thm2} $X$ is in the Fujiki class $\mathcal{C}$, hence it is a sGG-manifold. This being an open condition (see \cite{PU}), up to shrinking $B$ we can suppose that for every $b\in B$ the manifold $X_{b}$ is sGG.

Moreover, as $(X,\phi)\in\overline{\mathcal{M}}_{Z}^{hk}$, there is a sequence $\{b_{n}\}$ of points of $B$ converging to $0$ over which the fiber $X_{n}$ is irreducible hyperk\"ahler, and we can even suppose that $NS(X_{n})=0$. 

Consider a sequence $\{\alpha_{n}\}$ given by $\alpha_{n}\in\mathcal{C}_{X_{n}}$ converging to $\alpha$, hence either $\alpha_{n}\in\mathcal{C}_{X_{n}}^{+}$ (for all $n$), or $-\alpha_{n}\in\mathcal{C}_{X_{n}}^{+}$ (for all $n$). As recalled before, we have $\mathcal{C}_{X_{n}}^{+}\subseteq\mathcal{E}_{X_{n}}$: it follows that either $\alpha_{n}\in\mathcal{E}_{X_{n}}$ (for all $n$) or $-\alpha_{n}\in\mathcal{E}_{X_{n}}$ (for all $n$). By Theorem 5.9 of \cite{PU} we then conclude that either $\alpha\in\mathcal{E}_{X}$ or $-\alpha\in\mathcal{E}_{X}$.\endproof

The connected component of $\mathcal{C}_{X}$ contained in $\mathcal{E}_{X}$ will be denoted $\mathcal{C}^{+}_{X}$ and called \textit{positive cone} of $X$, in analogy with the hyperk\"ahler case. The other connected component of $\mathcal{C}_{X}$ will be denoted $\mathcal{C}^{-}_{X}$.

\subsection{Deformations and K\"ahler classes}

The first result we prove is the following:

\begin{prop}
\label{prop:defok}
Let $X$ be a compact, connected holomorphic symplectic manifold in the Fujiki class $\mathcal{C}$, which is limit of irreducible hyperk\"ahler manifolds. Let $\alpha\in\mathcal{C}_{X}$.
\begin{enumerate}
 \item If for every $\beta\in H^{2}(X,\mathbb{Z})$ we have $q(\alpha,\beta)\neq 0$, then there is $t\in T(\alpha)$ such that $\mathcal{X}_{t}$ is K\"ahler and either $\alpha_{t}$ or $-\alpha_{t}$ is a K\"ahler class on $\mathcal{X}_{t}$.
 \item If moreover $\alpha\in\mathcal{C}^{+}_{X}$, then $\alpha_{t}$ is K\"ahler.
\end{enumerate}
\end{prop}

\proof Let $Z$ be an irreducible hyperk\"ahler manifold which is deformation equivalent to $X$, and let $(\Lambda,q)$ be its Beauville lattice. Fix a marking $\phi$ on $X$ (which exists by Theorem \ref{thm:thm1}), and consider the point $(X,\phi)\in\mathcal{M}_{Z}$. As $X$ is limit of irreducible hyperk\"ahler manifolds, we have $(X,\phi)\in\overline{\mathcal{M}}^{hk}_{Z}$.

We first show that for the generic $t\in T(\alpha)$ the fiber $\mathcal{X}_{t}$ is in $\overline{\mathcal{M}}^{hk}_{Z}$. Let $\mathcal{X}\longrightarrow B$ be the Kuranishi family of $X$, and let $0$ be the point of $B$ over which the fiber $X_{0}$ is $X$. 

As $(X,\phi)\in\overline{\mathcal{M}}^{hk}_{Z}$, there is a sequence $\{b_{n}\}$ of points of $B$ verifying the following properties:
\begin{enumerate}
 \item the sequence $b_{n}$ converges to $0$ in $B$;
 \item for each $n$ the fiber $X_{n}$ of $\mathcal{X}$ over $b_{n}$ is an irreducible hyperk\"ahler manifold such that $H^{1,1}(X_{n})\cap H^{2}(X_{n},\mathbb{Z})=0$;
 \item for each $n$ there is $\alpha_{n}\in\mathcal{C}_{X_{n}}$ such that the sequence $\alpha_{n}$ converges to $\alpha$. 
\end{enumerate}
As $H^{1,1}(X_{n})\cap H^{2}(X_{n},\mathbb{Z})=0$, up to changing the sign of $\alpha$, and hence of $\alpha_{n}$, we can suppose that $\alpha_{n}\in\mathcal{K}_{X_{n}}$ for every $n$. We let $T_{n}$ be the twistor line of $\alpha_{n}$, which is a rational curve in $B$ passing through the point $(X_{n},\phi_{n})$.

As $(X_{n},\phi_{n})$ converges to $(X,\phi)$, and as $\alpha_{n}$ converges to $\alpha$, we see that the twistor lines $T_{n}$ converge to $T(\alpha)$. This means that if $t\in T(\alpha)$, there is a sequence $\{s_{t,n}\}$ of points of $B$ such that
\begin{enumerate}
 \item the sequence $\{s_{t,n}\}$ converges to $t$;
 \item for each $n$ we have $s_{t,n}\in T_{n}$.
\end{enumerate}
As $s_{t,n}\in T_{n}$, and as $T_{n}$ is the twistor line of the K\"ahler class $\alpha_{n}$, we see that the fiber $X_{s_{t,n}}$ of the twistor family of $\alpha_{n}$ over $s_{t,n}$ is an irreducible hyperk\"ahler manifold. As $s_{t,n}$ converges to $t$, we then see that $X_{t}$ is limit of irreducible hyperk\"ahler manifolds. This means that $(X_{t},\phi_{t})\in\overline{\mathcal{M}}^{hk}_{Z}$.

In particular, by Proposition \ref{prop:thm2equiv} this implies that $X_{t}$ is bimeromorphic to an irreducible hyperk\"ahler manifold $Y_{t}$. Now, by hypothesis we have $q_{X}(\alpha,\beta)\neq 0$ for each $\beta\in H^{2}(X,\mathbb{Z})$. This implies that $T(\alpha)$ does not intersect in 0 (and hence generically) any of the hypersurfaces $S_{\beta}$, i. e. for a generic $t\in T(\alpha)$ the period of $(X_{t},\phi_{t})$ is generic in $\Omega_{\Lambda}$.

As the periods of $(X_{t},\phi_{t})$ and $(Y_{t},\psi_{t})$ are equal, it follows that for a generic $t\in T(\alpha)$ the irreducible hyperk\"ahler manifold $Y_{t}$ is such that $H^{1,1}(Y_{t})\cap H^{2}(Y_{t},\mathbb{Z})=0$, so that $X_{t}$ and $Y_{t}$ are biholomorphic.

It follows that $X_{t}$ is irreducible hyperk\"ahler, and that $\mathcal{K}_{X_{t}}$ is one of the components of $\mathcal{C}_{X_{t}}$. As $\alpha_{t}\in\mathcal{C}_{X_{t}}$, it follows that either $\alpha_{t}$ or $-\alpha_{t}$ is a K\"ahler class on $X_{t}$.

Let us now suppose moreover that the class $\alpha$ is even pseudo-effective, and by contradiction that $\alpha_{t}$ is not a K\"ahler class. By what we just proved, it follows that $-\alpha_{t}$ is K\"ahler for generic $t\in T(\alpha)$. As $\mathcal{K}_{X_{t}}$ is contained in $\mathcal{E}_{X_{t}}$, we then have a family $-\alpha_{t}$ of pseudo-effective classes converging to $-\alpha$.

Now, by Theorem 5.9 of \cite{PU} (which we can apply as by the previous part of the proof the family $\mathcal{X}(\alpha)\longrightarrow T(\alpha)$ is a family of class $\mathcal{C}$ manifolds, and hence of sGG manifolds) we know that a limit of pseudo-effective classes along the family $\mathcal{X}(\alpha)$ is a pseudo-effective class on $X$.

This means that $-\alpha$ is a pseudo-effective class on $X$. As by hypothesis $\alpha$ is pseudo-effective too, it follows that $\alpha=0$,  which is not possible as $q_{X}(\alpha)>0$. This shows that is $\alpha$ is a positive pseudo-effective class such that $q_{X}(\alpha,\beta)\neq 0$ for every $\beta\in H^{2}(X,\mathbb{Z})$, then for a generic $t\in T(\alpha)$ the class $\alpha_{t}$ is K\"ahler.\endproof

We now use the previous Proposition to show the following, which is an improved version of Proposition \ref{prop:thm2equiv}.

\begin{prop}
\label{prop:bimpos}Let $X$ be a compact holomorphic symplectic manifold in the Fujiki class $\mathcal{C}$, which is limit of irreducible hyperk\"ahler manifolds, and let $\alpha\in\mathcal{C}^{+}_{X}$ be such that $q_{X}(\alpha,\beta)\neq 0$ for every $\beta\in H^{2}(X,\mathbb{Z})$. Then there is an irreducible hyperk\"ahler manifold $Y$ and a cycle $\Gamma=Z+\sum_{i}D_{i}$ in $X\times Y$ such that the following properties are verified:
\begin{enumerate}
 \item the cycle $Z$ defines a bimeromorphic map between $X$ end $Y$;
 \item the projections $D_{i}\longrightarrow X$ and $D_{i}\longrightarrow Y$ have positive dimensional fibers;
 \item the cycle $\Gamma$ defines a Hodge isometry $[\Gamma]_{*}$ between $H^{2}(X,\mathbb{Z})$ and $H^{2}(Y,\mathbb{Z})$;
 \item the class $[\Gamma]_{*}\alpha$ is K\"ahler.
\end{enumerate}
\end{prop}

\proof Consider the family $\kappa_{\alpha}:\mathcal{X}(\alpha)\longrightarrow T(\alpha)$. By Proposition \ref{prop:defok}, we know that for a generic $t\in T(\alpha)$ the fiber $X_{t}$ of $\kappa_{\alpha}$ over $t$ is an irreducible hyperk\"ahler manifold, and that $\alpha_{t}$ is a K\"ahler class on it.

Let $\mathcal{X}'\longrightarrow T(\alpha_{t})$ be the twistor family of $(X_{t},\alpha_{t})$, and notice that $\pi(T(\alpha))$ is identified with an open subset of $\pi(T(\alpha_{t}))$, and that for every $s\in T(\alpha_{t})$ the fiber $X'_{s}$ of $\mathcal{X}'$ over $s$ is K\"ahler.

Restricting the twistor family $\mathcal{X}'$ to such an open subset, we then find two families $\mathcal{X}(\alpha)\longrightarrow C$ and $\mathcal{X}'\longrightarrow C$ over the same base curve, which have isomorphic fibers over $t$, and the fibers of $\mathcal{X}'$ are all K\"ahler. We let $0\in C$ be the point over which the fiber of $\mathcal{X}$ is $X$, and we let $X'$ be the fiber of $\mathcal{X}'$ over $0$.

Both families are endowed with natural markings $\phi_{s}$ and $\phi'_{s}$ for each $s$, such that $(\phi'_{t})^{-1}\circ\phi_{t}$ is induced by the isomorphism $X_{t}\simeq X'_{t}$. The class $\alpha'_{s}:=(\phi'_{s})^{-1}\circ\phi_{s}(\alpha_{s})$ is a K\"ahler class on $X'_{s}$ for every $s\in C$. In particular the class $\alpha':=(\phi'_{0})^{-1}\circ\phi_{0}(\alpha)$ is K\"ahler on $X'$. 

By Lemma \ref{lem:border} the points $(X,\phi_{0})$ and $(X',\phi'_{0})$ are non-separated points in $\mathcal{M}_{Z}$. As $(X,\phi_{0})\in\overline{\mathcal{M}}^{hk}_{Z}$ and $(X',\phi'_{0})\in\mathcal{M}^{hk}_{Z}$, we can apply Proposition \ref{prop:thm2equiv} to show that there is a cycle $\Gamma=Z+\sum_{i}Y_{i}$ on $X\times X'$ such that
\begin{enumerate}
 \item the cycle $Z$ defines a bimeromorphic map between $X$ end $X'$;
 \item the projections $Y_{i}\longrightarrow X$ and $Y_{i}\longrightarrow X'$ have positive dimensional fibers;
 \item the cycle $\Gamma$ defines a Hodge isometry $[\Gamma]_{*}$ between $H^{2}(X,\mathbb{Z})$ and $H^{2}(X',\mathbb{Z})$;
\end{enumerate}
Notice that $[\Gamma]_{*}\alpha=\alpha'$, which is K\"ahler.\endproof 

\subsection{The proof of Theorem \ref{thm:thm3}}

We are now ready to prove Theorem \ref{thm:thm3}, namely that if $X$ is a compact, connected holomorphically symplectic $b_{2}-$manifold which is limit of irreducible hyperk\"ahler manifolds, any very general class $\alpha\in\mathcal{C}_{X}^{+}$ (i. e. $q_{X}(\alpha,\beta)\neq 0$ for every $\beta\in H^{2}(X,\mathbb{Z})$) such that $\alpha\cdot C>0$ for every rational curve $C$ on $X$ is a K\"ahler class on $X$, and in particular $X$ is K\"ahler.

\proof By Proposition \ref{prop:bimpos}, as $\alpha\in\mathcal{C}^{+}_{X}$ is such that $q_{X}(\alpha,\beta)\neq 0$ for every $\beta\in H^{2}(X,\mathbb{Z})$, then there is an irreducible hyperk\"ahler manifold $Y$ and a cycle $\Gamma=Z+\sum_{i}D_{i}$ in $X\times Y$ such that the following properties are verified:
\begin{enumerate}
 \item the cycle $Z$ defines a bimeromorphic map between $X$ end $Y$;
 \item the projections $D_{i}\longrightarrow X$ and $D_{i}\longrightarrow Y$ have positive dimensional fibers;
 \item the cycle $\Gamma$ defines a Hodge isometry $[\Gamma]_{*}$ between $H^{2}(X,\mathbb{Z})$ and $H^{2}(Y,\mathbb{Z})$;
 \item the class $\alpha':=[\Gamma]_{*}\alpha$ is K\"ahler on $Y$.
\end{enumerate}

The argument used in the proof of Theorem 2.5 of \cite{Hu2} shows that since $[\Gamma]_{*}\alpha$ is a K\"ahler class on $Y$ and $\alpha\cdot C>0$ for every rational curve $C$ on $X$, then all the irreducible components $D_{i}$ of $\Gamma$ which are contracted by the projection $p_{X}$ of $X\times Y$ to $X$ are such that the codimension in $X$ of $p_{X}(D_{i})$ is at least 2. By Lemma 2.2 of \cite{Hu99} it then follows that the morphisms $[D_{i}]_{*}:H^{2}(Y,\mathbb{Z})\longrightarrow H^{2}(X,\mathbb{Z})$ are all trivial. As a consequence, we have $\alpha=[\Gamma]_{*}\alpha'=[Z]_{*}\alpha'$. 

We let $f:Y\dashrightarrow X$ be a bimeromorphism whose graph is $Z$. As $\alpha'$ is K\"ahler, then for every rational curve $C'$ in $Y$ we have $$\int_{C'}\alpha'>0.$$Notice that $\alpha'=f^{*}\alpha$, so that we have $$\int_{C}\alpha>0,\,\,\,\,\,\,\,\,\int_{C'}f^{*}\alpha>0$$for every rational curve $C$ in $X$ and every rational curve $C'$ in $Y$. By Proposition 2.1 of \cite{Hu2} it follows that $f$ extends to an isomorphism, and $\alpha$ is then a K\"ahler class.\endproof

\section{K\"ahlerness of moduli spaces of sheaves}

This last section is devoted to the proof of Theorem \ref{thm:main}. Hence, we let $S$ be a K3 surface, $v\in H^{2*}(S,\mathbb{Z})$ be of the form $v=(r,\xi,a)$ where $r>0$ and $\xi\in NS(S)$ are prime to each other. Moreover, we let $\omega$ be a K\"ahler class on $S$, which we suppose to be $v-$generic.

We want to show that if the moduli space $M:=M^{\mu}_{v}(S,\omega)$ is a $b_{2}-$manifold, then it is K\"ahler. To do so, we apply Theorem  \ref{thm:thm3} to $M$: we then need to prove that $M$ is a compact, connected holomorphically symplectic $b_{2}-$manifold which is limit of irreducible hyperk\"ahler manifolds, and we need to provide a very generic class $\alpha\in\mathcal{C}^{+}_{M}$ such that $\alpha\cdot C>0$ for every rational curve $C$ in $M$.

We will always assume that $v^{2}\geq 2$, as the cases $v^{2}\leq 0$ are already known: if $v^{2}<-2$ we have $M=\emptyset$; if $v^{2}=-2$ then $M$ is a point; if $v^{2}=0$, then $M$ is a K3 surface by Corollary 5.3 of \cite{PT}.

\paragraph{A.   The moduli space $M$ is a compact, connected holomorphically symplectic $b_{2}-$manifold which is limit of irreducible hyperk\"ahler manifolds.} The fact that $M$ is a compact, holomorphically symplectic manifold is due to Toma (see Remark 4.5 of \cite{T}). The connectedness is given by Proposition 4.24 of \cite{PT}. The fact that $M$ is a $b_{2}-$manifold is supposed to hold true.

We are then left to prove that $M$ is a limit of irreducible hyperk\"ahler manifolds. To do so, let $\mathcal{S}\longrightarrow B$ be the Kuranishi family of the K3 surface $S$, where $B$ is a complex manifold of dimension 20. Let $B_{\xi}\subseteq B$ be the subvariety of $B$ given by those $b\in B$ such that $\xi\in NS(S_{b})$. Similarily, let $B_{\omega}\subseteq B$ be the subvariety of $B$ given by those $b\in B$ such that the class $\omega\in H^{1,1}(S_{b},\mathbb{R})$. Moreover, we let $B_{\xi,\omega}:=B_{\xi}\cap B_{\omega}$.

Recall that $B_{\xi}$ and $B_{\omega}$ are smooth hypersurfaces of $B$, and as $\xi$ and $\omega$ are linearily independent, then $B_{\xi}$ and $B_{\omega}$ intersect transversally, so that $B_{\xi,\omega}$ is a smooth analytic subset of $B$ (of positive dimension). By Theorem 3.5 of \cite{Hu99}, the subset $B^{p}_{\xi,\omega}$ of $B_{\xi,\omega}$ given by those $b$ such that $S_{b}$ is projective is dense in $B_{\xi,\omega}$.

We now consider the restriction $\mathcal{S}':=\mathcal{S}_{|B_{\xi,\omega}}$, together with a morphism $\mathcal{S}'\longrightarrow B_{\xi,\omega}$. We suppose $0\in B_{\xi,\omega}$ be such that $S_{0}\sim S$. 

Recall that $\omega$ is a K\"ahler class on $S$: as K\"ahlerness is an open property in the analytic topology, there is an analytic open subset $D'\subset B_{\xi,\omega}$ containing $0$ and such that for every $b\in D'$ the class $\omega$ is K\"ahler on $S_{b}$. We then can consider the relative moduli space $\mathcal{M}\longrightarrow D'$, whose fiber over $b$ is the moduli space $M_{b}=M^{\mu}_{v}(S_{b},\omega)$ of $\mu_{\omega}-$stable sheaves on $S_{b}$ whose Mukai vector is $v$.

As the $v-$genericity is an open property in the analytic topology, there is an open subset $D$ of $D'$ such that for every $b\in D$ the class $\omega$ is $v-$generic. We then consider the restriction $\mathcal{S}_{D}$ of $\mathcal{S}$ to $D$, together with a morphism $\mathcal{S}_{D}\longrightarrow D$. For every $d\in D$ the K3 surface $S_{d}$ comes equipped with a Mukai vector $v=(r,\xi,a)$ and a $v-$generic polarization $\omega$.

As a consequence, the restriction of $\mathcal{M}$ to $D$, denoted $\mathcal{M}_{D}$ is such that for every $d\in D$ the fiber $M_{d}$ is a compact, connected complex manifold (see again Proposition 4.24 of \cite{PT}). The morphism $\mathcal{M}_{D}\longrightarrow D$ being submersive, it follows that the family $\mathcal{M}_{D}\longrightarrow D$ is a smooth and proper family, whose fiber over $0$ is $M^{\mu}_{v}(S,\omega)$. 

Now, as $B^{p}_{\xi,\omega}$ is dense in $B_{\xi,\omega}$, it follows that $B^{p}_{\xi,\omega}\cap D$ is dense in $D$. Hence, for the generic point $d\in D$ the fiber $S_{d}$ is a projective K3 surface, and the fiber $M_{d}$ is an irreducible hyperk\"ahler manifold (see Theorem 3.4 of \cite{PT}). Hence $M^{\mu}_{v}(S,\omega)$ is limit of irreducible hyperk\"ahler manifolds.

\begin{oss}
\label{oss:fujiki}
{\rm As a consequence of what we just proved, by Theorem \ref{thm:thm2} we see that if $M^{\mu}_{v}(S,\omega)$ verifies the assumptions of Theorem \ref{thm:main}, then it is in the Fujiki class $\mathcal{C}$.}
\end{oss}

\paragraph{B.   A very generic class $\alpha\in\mathcal{C}^{+}_{M}$ such that $\alpha\cdot C>0$ for every rational curve $C$ of $M$.} In order to prove Theorem \ref{thm:main}, by the previous part of this section we just need to find such a class $\alpha$ on $M$. To do so, recall that by \cite{PT} there is a morphism $$\lambda_{v}:v^{\perp}\longrightarrow H^{2}(M,\mathbb{Z})$$which is an isometry (since $v^{2}\geq 2$) with respect to the Mukai pairing on $v^{\perp}$ and the Beauville form of $M$.

Moreover, by the previous paragraph the moduli space $M$ is in the Fujiki class $\mathcal{C}$, hence $H^{2}(M,\mathbb{Z})$ has a Hodge decomposition, and $\lambda_{v}$ is a Hodge morphism. In particular, $\lambda_{v}$ is a Hodge isometry. This remains true if we tensor with $\mathbb{R}$, and we get a Hodge isometry $$\lambda_{v}:v^{\perp}\otimes\mathbb{R}\longrightarrow H^{2}(M,\mathbb{R}).$$We will then construct the desired class $\alpha$ by taking an appropriate element of $v^{\perp}\otimes\mathbb{R}$.

The choice we make is the following: let $m\in\mathbb{N}$ and $$\alpha_{m,\omega}:=(-r,-mr\omega,a+m\omega\cdot\xi).$$First of all, we remark that $\alpha_{m,\omega}\in v^{\perp}\otimes\mathbb{R}$, as $$(v,\alpha_{m,\omega})_{S}=mr\omega\cdot\xi-r(a+m\omega\cdot\xi)+ra=0.$$Moreover, as $\omega$ is a real $(1,1)-$class on $S$, the class $\alpha_{m,\omega}$ is a real $(1,1)-$class orthogonal to $v$.

It then follows that $$\alpha:=\lambda_{v}(\alpha_{m,\omega})\in H^{1,1}(M,\mathbb{R}).$$We prove that $\alpha$ is a very general class in $\mathcal{C}_{M}^{+}$ such that $\alpha\cdot C>0$ for every rational curve $C$ on $M$.

We start by showing that $\alpha$ is very general in $H^{1,1}(M,\mathbb{R})$

\begin{lem}
\label{lem:verygeneric} 
If $\omega$ is sufficiently generic, then for every $\beta\in H^{2}(M,\mathbb{Z})$ we have $q_{M}(\alpha,\beta)\neq 0$. 
\end{lem}

\proof As $\beta\in H^{2}(M,\mathbb{Z})$, there are $s,b\in\mathbb{Z}$ and $D\in H^{2}(S,\mathbb{Z})$ such that $D\cdot\xi=sa+rb$ (i. e. $\gamma:=(s,D,b)\in v^{\perp}$) and $\beta=\lambda_{v}(\gamma)$. 

It follows that $$q_{M}(\alpha,\beta)=q_{M}(\lambda_{v}(\alpha_{m,\omega}),\lambda_{v}(\gamma))=(\alpha_{m,\omega},\gamma)_{S}=$$ $$=-m\omega\cdot(rD+s\xi)+rb-sa.$$Suppose that $q(\alpha,\beta)=0$: this is then equivalent to $$\omega\cdot(rD+s\xi)=\frac{D\cdot\xi-2sa}{m},$$which means that $\omega$ is on some hyperplane in $H^{2}(S,\mathbb{R})$ associated to $D$. As the family of these hyperplanes is countable (since the family of $D\in H^{2}(S,\mathbb{Z})$ is countable), and as $\omega$ is sufficiently generic, we see that $q(\alpha,\beta)\neq 0$ for every $\beta\in NS(M)$.\endproof

We notice that we can move $\omega$ in the $v-$chamber of the K\"ahler cone of $S$ where it lies without changing $M$ (see Proposition 3.2 of \cite{PT}), hence we can always suppose that $\omega$ is sufficiently generic, and hence that $\alpha$ is very generic.
 
We now show that $q_{M}(\alpha)>0$ and that it is a pseudo-effective class on $M$.

\begin{lem}
\label{lem:positivecone}
If $m\gg 0$ we have $\alpha\in\mathcal{C}_{M}^{+}$. 
\end{lem}

\proof We first prove that $\alpha\in\mathcal{C}_{M}$, and then that $\alpha\in\mathcal{E}_{M}$.

We have $$q_{M}(\alpha)=q_{M}(\lambda_{v}(\alpha_{m,\omega}))=(\alpha_{m,\omega},\alpha_{m,\omega})_{S}=$$ $$=m^{2}r^{2}\omega^{2}+2ra-2mr\omega\cdot\xi.$$As $m\gg 0$ and $\omega^{2}>0$ (since $\omega$ is K\"ahler on $S$), we then see that $q_{M}(\alpha)>0$, i. e. $\alpha\in\mathcal{C}_{M}$.

We now have to show that $\alpha\in\mathcal{E}_{M}$. To show this, consider the deformation $$\mathcal{M}\longrightarrow B_{\xi,\omega}$$we introduced in the previous paragraph. We let $0\in B_{\xi,\omega}$ be the point over which the fiber is $M^{\mu}_{v}(S,\omega)$. For a generic $b\in B_{\xi,\omega}$ the fiber is $M_{b}=M^{\mu}_{v}(S_{b},\omega)$ where $S_{b}$ is a projective K3 surface, so that the fiber is a projective irreducible hyperk\"ahler manifold.

Notice that $\omega$ is still a $v-$generic K\"ahler class on $S_{b}$, and the class $\alpha$ is still in $\mathcal{C}_{M_{b}}$, and this for every $b\in B_{\xi,\omega}$. We write $\alpha_{0}:=\alpha$

Now, as shown in Remark 3.5 of \cite{PT}, in the same $v-$chamber where $\omega$ lies there is a class of the form $\omega'=c_{1}(H)$ for some ample line bundle $H$ on $S_{b}$. We let $\alpha_{1}:=\lambda_{v}(\alpha_{m,\omega'})$. Moreover, for every $t\in[0,1]$ we let $\omega_{t}:=(1-t)\omega+t\omega'$, which is a segment contained in the $v-$chamber where $\omega$ and $\omega'$ are, and we let $\alpha_{t}:=\lambda_{v}(\alpha_{m,\omega_{t}})$. By linearity of $\lambda_{v}$, we have $$\alpha_{t}=(1-t)\alpha_{0}+t\alpha_{1},$$and the image of the map $\alpha:[0,1]\longrightarrow H^{1,1}(M_{b},\mathbb{R})$ defined by letting $\alpha(t):=\alpha_{t}$ is a segment in $\mathcal{C}_{M_{b}}$.

Our aim is to show that $\alpha\in \mathcal{E}_{M}$. As the family $\mathcal{M}\longrightarrow B$ is a family of manifolds in the Fujiki class $\mathcal{C}$ by the previous paragraph, by Theorem 5.9 of \cite{PU} it is sufficient to show that $\alpha_{0}\in\mathcal{E}_{M_{b}}$ for a generic $b$ around $0$. As for the generic $b$ around 0 we have that $M_{b}$ is irreducible hyperk\"ahler, this is equivalent to show that $\alpha_{0}\in\mathcal{C}_{M_{b}}^{+}$.

As $\mathcal{C}^{+}_{M_{b}}$ is a convex cone and the segment $[\alpha_{0},\alpha_{1}]$ is contained in $\mathcal{C}_{M_{b}}$, to show that $\alpha_{0}\in\mathcal{C}^{+}_{M_{b}}$ it is sufficient to show that $\alpha_{1}\in\mathcal{C}_{M_{b}}^{+}$.

But now, as $S_{b}$ is projective we can use a general construction presented in \cite{HL}: if $H$ is a $v-$generic ample line bundle on $S_{b}$, we can construct an ample line bundle $L(H)$ on $M^{\mu}_{v}(S_{b},H)$, and we have $c_{1}(L(H))=\lambda_{v}(\alpha_{m,c_{1}(H)})$. As $M^{\mu}_{v}(S_{b},\omega)=M^{\mu}_{v}(S_{b},H)$ (as $\omega$ and $c_{1}(H)$ are in the same $v-$chamber, by \cite{PT}), the class $\lambda_{v}(\alpha_{m,\omega'})$ is an ample class on $M^{\mu}_{v}(S_{b},\omega)$. It then lies in the K\"ahler cone of $M_{b}$, and hence in $\mathcal{C}^{+}_{M_{b}}$.\endproof

In conclusion, we have shown that up to choosing $m\gg 0$ and $\omega$ sufficiently generic, the class $\alpha$ is a very generic class in $\mathcal{C}_{M}^{+}$. We are left to show that $\alpha\cdot C>0$ for every rational curve $C$ in $M$.

\begin{lem}
\label{lem:positivecurve}
If $m\gg 0$ and $\omega$ is sufficiently generic, we have $\alpha\cdot C>0$ for every rational curve $C$ on $M$.
\end{lem}

\proof Let $[C]\in H^{2n-1,2n-1}(M,\mathbb{Z})$, and let $\beta_{C}\in NS(M)$ be the dual of $[C]$, so that $$\alpha\cdot C=q_{M}(\alpha,\beta_{C}).$$We then just need to prove that $q_{M}(\alpha,\beta_{C})>0$ for every rational curve $C$ on $M$.

Let $\mathcal{S}\longrightarrow B$ be the Kuranishi family of $S$, and let $0\in B$ be such that $S_{0}=S$. We let $B_{C}$ be the subset of $B$ of those $b\in B$ such that $\beta_{C}\in NS(S_{b})$, i. e. $C$ is a rational curve on $S_{b}$. Consider the intersection $D_{C}:=D\cap B_{C}$, which is an analytic subset of $D$, whose generic point $d$ is such that $S_{d}$ is a projective K3 surface.

We let $\mathcal{M}_{C}$ be the restriction of the relative moduli space $\mathcal{M}\longrightarrow D$ to $D_{C}$. Notice that for every $d\in D_{C}$ we have the class $\alpha\in\mathcal{C}_{M_{d}}$ and the rational curve $C$ on $M_{d}$. As the intersection product of $\alpha$ with $C$ is constant along $D_{C}$, it is sufficient to show that $$q_{M_{d}}(\alpha,\beta_{C})>0$$for some $d\in D_{C}$.

As $\beta_{C}\in NS(M_{d})$, there are $s,b\in\mathbb{Z}$ and $\zeta\in NS(S_{d})$ such that $\gamma:=(s,\zeta,b)\in v^{\perp}$ and $\beta_{C}=\lambda_{v}(\gamma)$. As $\lambda_{v}$ is an isometry, we have $$q_{M_{d}}(\alpha,\beta_{C})=q_{M_{d}}(\lambda_{v}(\alpha_{m,\omega}),\lambda_{v}(\gamma))=(\alpha_{m,\omega},\gamma)_{S_{d}}.$$It is then sufficient to show that $(\alpha_{m,\omega},\gamma)_{S_{d}}>0$.

Now, by Lemma 3.3 of \cite{PT} there is an ample class $\omega'$ on $S_{d}$ which is in the same $v-$chamber of $\omega$ and such that for every $\eta\in NS(S_{d})$ we have $\omega\cdot\eta=\omega'\cdot\eta$. Then we have $\alpha_{m,\omega'}\in v^{\perp}$, and $$(\alpha_{m,\omega},\gamma)_{S_{d}}=(\alpha_{m,\omega'},\gamma)_{S_{d}}.$$It is then sufficient to show that $(\alpha_{m,\omega'},\gamma)_{S_{d}}>0$.

To do so, consider a rational $\omega''$ in a neighborhood of $\omega'$ in the ample cone of $S_{d}$, and let $p\in\mathbb{N}$ and $H$ an ample line bundle on $S_{d}$ such that $p\omega''=c_{1}(H)$. As we can choose $m\gg 0$, we can suppose that $m=m'p$ for some very big $m'\in\mathbb{N}$. As $H$ is $v-$generic we have that $\lambda_{v}(\alpha_{m',c_{1}(H)})$ is the first Chern class of an ample line bundle, so that $$\lambda_{v}(\alpha_{m',c_{1}(H)})\cdot C>0.$$It follows that $$(\alpha_{m',c_{1}(H)},\gamma)_{S_{d}}>0.$$

Now, notice that $$\alpha_{m',c_{1}(H)}=(-r,m'rc_{1}(H),a+m'c_{1}(H)\cdot\xi)=$$ $$=(-r,m'prc_{1}(H)/p,a+m'p\xi\cdot c_{1}(H)/p)=$$ $$=(-r,mr\omega'',a+m\omega''\cdot\xi)=\alpha_{m,\omega''}.$$Hence we get $$(\alpha_{m,\omega''},\gamma)_{S_{d}}>0$$As this is true for all rational classes $\omega''$ in a neighborhood of $\omega'$, this implies that $$(\alpha_{m,\omega'},\gamma)_{S_{d}}\geq 0.$$

As we saw before, this implies that $\alpha\cdot C\geq 0$ for every rational curve $C$ in $M$. But as $\beta_{C}\in NS(M)$ and $\alpha\cdot C=q_{M}(\alpha,\beta_{C})$, and as we know that $\alpha$ is very generic by Lemma \ref{lem:verygeneric}, it follows that $\alpha\cdot C\neq 0$. In conclusion we have $\alpha\cdot C>0$, and we are done.\endproof

Now, by paragraph A we know that if $M^{\mu}_{v}(S,\omega)$ verifies the assumptions of Theorem \ref{thm:main}, then it is a compact, connected holomorphically symplectic $b_{2}-$manifold which is limit of irreducible hyperk\"ahler manifolds. Moreover, Lemmas \ref{lem:verygeneric}, \ref{lem:positivecone} and \ref{lem:positivecurve} we know that the class $\alpha$ is a very generic class in $\mathcal{C}_{M}^{+}$ such that $\alpha\cdot C>0$ for every rational curve $C$ of $M^{\mu}_{v}(S,\omega)$. 

We can then apply Theorem \ref{thm:thm3} to conclude that $M^{\mu}_{v}(S,\omega)$ is a K\"ahler manifold (actually irreducible hyperk\"ahler), and that $\alpha$ is a K\"ahler class on it. This concludes the proof of Theorem \ref{thm:main}.

\par\bigskip
\par\bigskip

Institut \'Elie Cartan, UMR 7502, Universit\'e de Lorraine, CNRS, INRIA, Boulevard des Aiguillettes, B.P. 70239, 54506 Vandoeuvre-l\`es-Nancy Cedex, France.

E-mail: arvid.perego@univ-lorraine.fr

\end{document}